\documentclass[11 pt]{amsart}
\usepackage{amsmath}
\usepackage{amsfonts}
\usepackage{amssymb}
\usepackage{color}
\theoremstyle{remark}

\theoremstyle{definition}

\newcommand{\average}{{\mathchoice {\kern1ex\vcenter{\hrule height.4pt
width 6pt
depth0pt} \kern-9.7pt} {\kern1ex\vcenter{\hrule height.4pt width 4.3pt
depth0pt}
\kern-7pt} {} {} }}

\def\vint_#1{\mathchoice%
          {\mathop{\kern 0.2em\vrule width 0.6em height 0.69678ex depth -0.58065ex
                  \kern -0.8em \intop}\nolimits_{\kern -0.4em#1}}%
          {\mathop{\kern 0.1em\vrule width 0.5em height 0.69678ex depth -0.60387ex
                  \kern -0.6em \intop}\nolimits_{#1}}%
          {\mathop{\kern 0.1em\vrule width 0.5em height 0.69678ex
              depth -0.60387ex
                  \kern -0.6em \intop}\nolimits_{#1}}%
          {\mathop{\kern 0.1em\vrule width 0.5em height 0.69678ex depth -0.60387ex
                  \kern -0.6em \intop}\nolimits_{#1}}}
\def\vintslides_#1{\mathchoice%
          {\mathop{\kern 0.1em\vrule width 0.5em height 0.697ex depth -0.581ex
                  \kern -0.6em \intop}\nolimits_{\kern -0.4em#1}}%
          {\mathop{\kern 0.1em\vrule width 0.3em height 0.697ex depth -0.604ex
                  \kern -0.4em \intop}\nolimits_{#1}}%
          {\mathop{\kern 0.1em\vrule width 0.3em height 0.697ex depth -0.604ex
                  \kern -0.4em \intop}\nolimits_{#1}}%
          {\mathop{\kern 0.1em\vrule width 0.3em height 0.697ex depth -0.604ex
                  \kern -0.4em \intop}\nolimits_{#1}}}

\newcommand{\kintint}[2]{\mathchoice%
          {\mathop{\kern 0.2em\vrule width 0.6em height 0.69678ex depth -0.58065ex
                  \kern -0.8em \intop}\nolimits_{\kern -0.45em#1}^{#2}}%
          {\mathop{\kern 0.1em\vrule width 0.5em height 0.69678ex depth -0.60387ex
                  \kern -0.6em \intop}\nolimits_{#1}^{#2}}%
          {\mathop{\kern 0.1em\vrule width 0.5em height 0.69678ex depth -0.60387ex
                  \kern -0.6em \intop}\nolimits_{#1}^{#2}}%
          {\mathop{\kern 0.1em\vrule width 0.5em height 0.69678ex depth -0.60387ex
                  \kern -0.6em \intop}\nolimits_{#1}^{#2}}}

\makeatletter
\def\cleardoublepage{\clearpage\if@twoside \ifodd\c@page\else
\hbox{}
\thispagestyle{empty}
\newpage
\if@twocolumn\hbox{}\newpage\fi\fi\fi}
\makeatother
\title{
Viscosity and variational approaches in free boundary problems: A Volume in Honor of Sandro Salsa}
\author{Daniela De Silva}
\address{Daniela De Silva: Department of Mathematics, Barnard College, Columbia University, New York, NY 10027, USA}
\email{desilva@math.columbia.edu }
\author{Fausto Ferrari}
\address{Fausto Ferrari: Dipartimento di Matematica\\ Universit\`a di Bologna\\ Piazza di Porta S.Donato 5\\ 40126, Bologna-Italy}
\email{fausto.ferrari@unibo.it }

\thanks{The authors wish to thank Donatella Danielli, Editor-in-Chief of La Matematica, for hosting this special volume and for her support throughout the many activities involved in its preparation and production.
 F. F. was partially supported by the Italian Ministry for Universities and Scientific and Technological Research MURST grant PRIN 2022 7HX33Z- CUP J53D23003610006, Pattern formation in nonlinear phenomena; Istituto Nazionale di Alta Matematica INdAM-GNAMPA 2024 project: Free boundary problems in noncommutative structures and degenerate operators CUP E53C23001670001; Dipartimento di Matematica dell'Universit\`a di Bologna, Italy.}
\date{\today}
\usepackage{graphicx}
\begin{document}
\maketitle

An eminent figure in modern mathematics, Sandro Salsa is Professor Emeritus at the ``Politecnico di Milano" (Milan, Italy). His career spans decades of research, teaching, and academic leadership that have left a lasting mark on the field of mathematical analysis, particularly through his extensive work on partial differential equations (PDEs) and free boundary problems.

Salsa graduated in Mathematics from the University of Milan in 1972 and joined the Politecnico in 1987 as a full professor of Mathematical Analysis, a position he held with distinction for over thirty years. As Director of the Department of Mathematics from 1999 to 2008, his leadership expanded the institution's international standing. Driven by an interdisciplinary outlook, he used his tenure to help establish two major initiatives: the MOX Laboratory for Modeling and Scientific Computing---which he helped develop into a highly regarded center---and the Mathematical Engineering program. This early focus on 
mathematical engineering anticipated, by more than a decade, the widespread modern interest in scientific computing and data science across industry and academia.

His impact as an educator is equally clear. Salsa has written widely used textbooks on mathematical analysis and PDEs that have become standard references in university education due to their exceptional clarity. Generations of students in engineering, mathematics, and economics have benefited from his clear teaching style, passion for the subject, and educational vision.

Alongside this institutional and educational legacy, Salsa’s research career has centered on nonlinear, degenerate, and singular elliptic and parabolic PDEs, making significant contributions to fluid dynamics, optimization, and free boundary problems. His work includes influential papers on regularity theory, the boundary behavior of elliptic operators, and obstacle problems, often co-authored with leading international colleagues and published in top-tier journals such as 
\textit{Acta Mathematica}, \textit{Annals of Mathematics}, and \textit{Communications 
on Pure and Applied Mathematics}. 

Among the many papers Salsa has written, several have introduced highly influential methods that anticipated developments in the field. In \cite{CFMS}, for instance, the authors extended the classical results of R. A. Hunt and R. L. Wheeden \cite{HW} on the boundary behavior of positive harmonic functions in Lipschitz domains to a class of second-order elliptic operators in divergence form, establishing a 
boundary Harnack principle. This work introduced a variety of fundamental 
mathematical tools that remain widely used today. From this perspective, it is important to highlight the foundational, though sometimes implicit, influence of Eugene Fabes and the Minneapolis school of PDEs during the 1980s on Sandro Salsa and his contemporaries. In this framework, it is worth to cite the contribution given in \cite{JK} about the fine properties of the harmonic functions in NTA domains. Within this network, Salsa developed highly productive scientific relationships, most significantly collaborating on the properties of the heat operator and broader classes of parabolic operators, as seen in papers such as \cite{FGS1, FGS2, FGMS}. 
\begin{figure}
\begin{center}
\includegraphics[width=5in, scale=0.5]{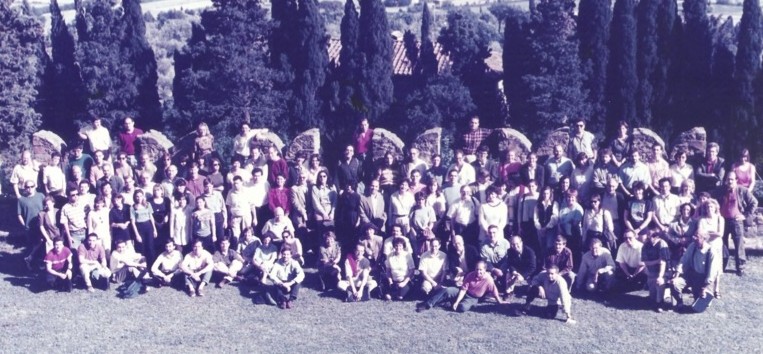}
\end{center}
\nopagebreak{
{\begin{center} Meeting in honor of Eugene Fabes, Cortona (Arezzo) 1996. Sandro Salsa is standing at the far right. \end{center}}}
\end{figure}

In \cite{FGS2}, the authors notably investigated the classical parabolic Harnack inequality for nonnegative solutions $u$ of 
$$
\text{div}(A(x,t)\nabla u)-\frac{\partial u}{\partial t}=0
$$
in $D_+ := D \times (0, +\infty)$, where $A$ is uniformly elliptic. While the standard inequality states that the values of $u$ at an earlier time $t_1$ are controlled by its values at a later time $t_2 > t_1$--intuitively interpreted as \textit{the minimum temperature in the future controlling the maximum temperature 
in the past}--the authors demonstrated that this relation can be reversed under specific boundary conditions. They showed that when a nonnegative solution vanishes on the entire lateral boundary $S_+ := \partial D \times (0, +\infty)$ of a bounded domain, the forward inequality can turn into a backward Harnack inequality, where the values of $u$ at a later time $t_2$ are bounded by its values at a earlier time $t_1$.

Salsa’s research on parabolic operators culminated in a prominent trilogy of papers, \cite{ACS1, ACS2, ACS3}, co-authored with I. Athanasopoulos and L. Caffarelli. This work established foundational regularity results for the free boundary in the Stefan problem, which remain benchmarks in the literature today. Building on these breakthroughs, Caffarelli and Salsa undertook a major effort to 
systematically organize the results of free boundary problems where viscosity theory plays a central role in the proofs. The resulting landmark monograph, \textit{A Geometric Approach to Free Boundary Problems} \cite{CS}, was published in 2005 by the American Mathematical Society. This work is widely recognized as a definitive milestone in the literature, leaving the next generation of mathematicians with an essential foundation on which to base and continue research in the field.
\begin{figure}
\begin{center}
\includegraphics[width=3in, scale=0.5]{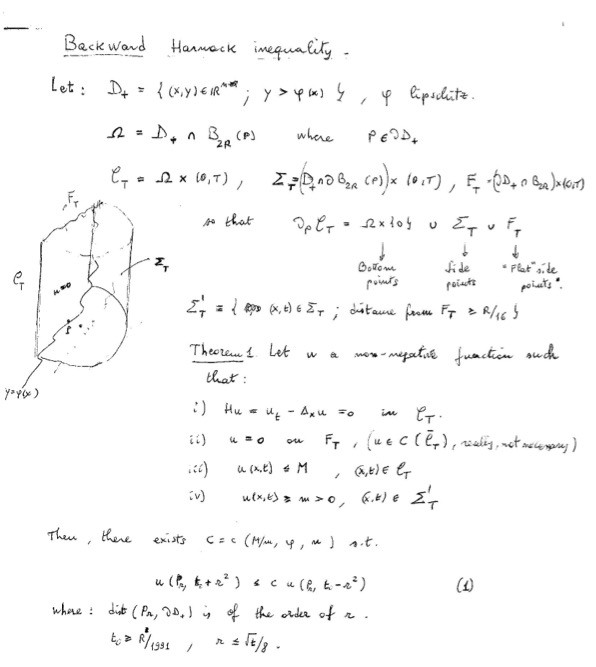}
\end{center}
\nopagebreak{
{\begin{center} Backward Harnack Inequality in Sandro Salsa's notes  \end{center}}}
\end{figure}

\subsection*{Free Boundary Problems: An Overview}

Free Boundary Problems (FBPs) are a class of differential equations where the domain of definition depends on the solution itself and must be determined as part of the problem. Arising across diverse physical contexts---such as flame propagation, fluid mechanics, optimal stopping in finance, and tumor growth---these problems form a rich interface between the applied sciences and pure mathematical analysis. 

Historically, the modern discipline traces its roots to the classic Stefan 
problem, which describes the phase-transition between a liquid and a solid, see for instance \cite{Friedman} and \cite{Kind-Stam-book} to deepen the subject or also \cite{CSV} for a short overview. Considered as early as 1831 by G. Lam\'e and B. Clapeyron in relation to ice formation, it was systematically analyzed by J. Stefan around 1890 using thermodynamic principles. A rigorous analytical treatment, however, required the modern framework of weak solutions and nonlinear analysis, culminating in the foundational existence and uniqueness results by S. Kamenomostkaya \cite{Kam} and O. Oleinik. 
 \begin{figure}
\begin{center}
\includegraphics[width=3in, scale=0.5]{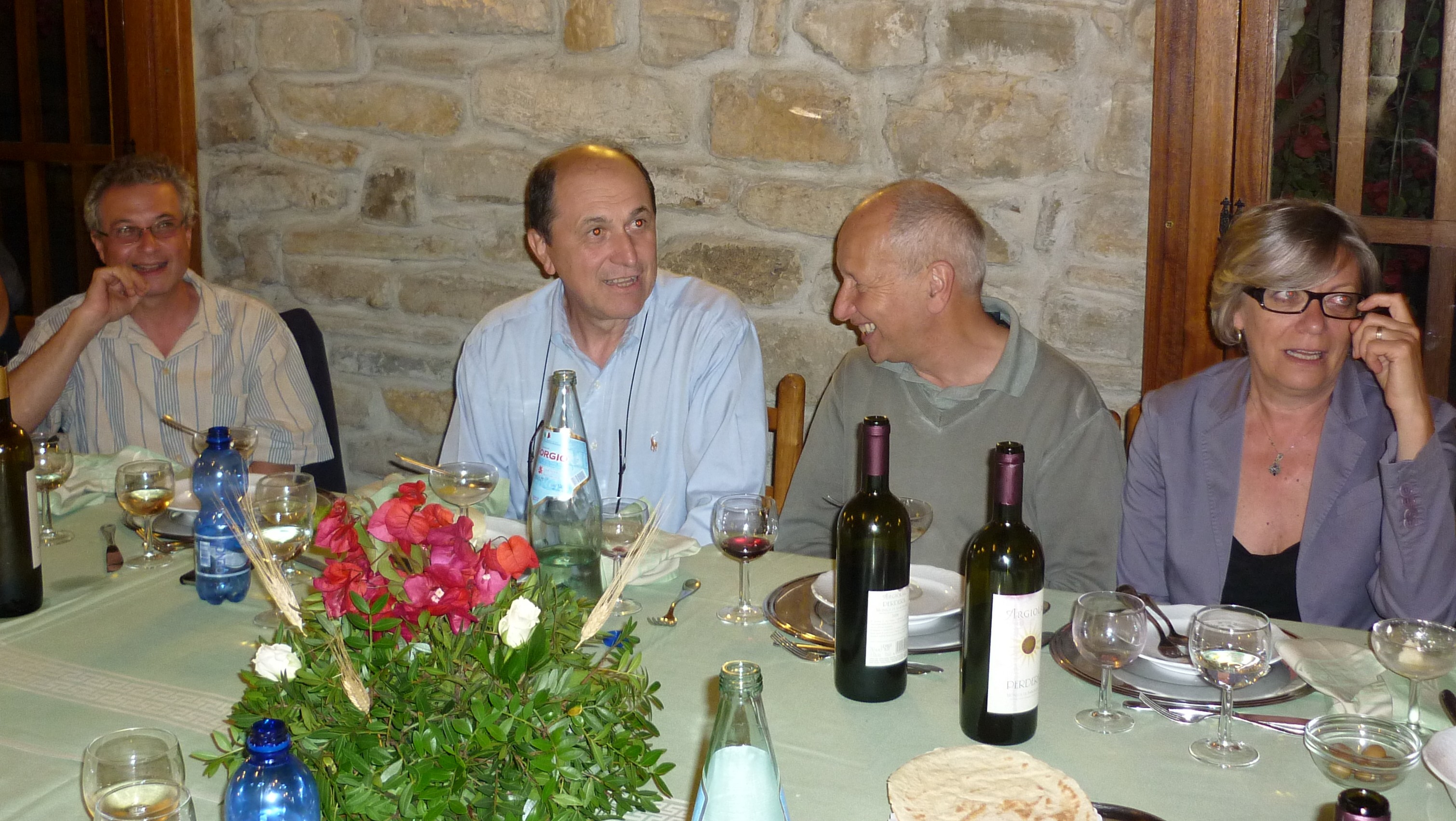}
\end{center}
\nopagebreak{
{\begin{center} From the left:  Carlos Kenig, Sandro Salsa, Luis Caffarelli, Anna Aloe (Salsa's wife). Neighborhood of Cagliari, 2010  \end{center}}}
\end{figure}

In its simplest, homogeneous, one-phase formulation, the Stefan problem can be expressed as:
\begin{equation}\label{eq:StefanNHSS}
\begin{cases}
  \partial_{t} u - \Delta u  =  0 & \text{in} \,\, \Omega^+(u):=\left(\Omega \times [0,T] \right) \cap \{ u > 0  \},\\
  \partial_{t} u = |\nabla u|^2 & \text{on} \,\, F(u):=\left(\Omega \times [0,T] \right) \cap \partial\{ u > 0  \},
\end{cases}
\end{equation}
where $\Omega \subset \mathbb{R}^{n}$ and $T > 0$. Here, $u$ denotes the temperature, and energy conservation is enforced at the interface. Geometrically, if $V_\nu := \frac{\partial_t u}{|\nabla u|}$ denotes the velocity of the free boundary $F(u)$ along the spatial normal $\nu$, the free boundary condition reads $V_\nu = |\nabla u|$, which means that the normal velocity of the interface is dictated by the temperature gradient.

A second classical archetype is the Bernoulli problem, which originates from Daniel Bernoulli's hydrodynamic studies on constant-pressure interfaces in steady, irrotational, and inviscid flows. Mathematically, it requires determining a harmonic potential $\psi$ (the stream function) inside an unknown domain where a gradient condition is prescribed on the free boundary $F(\psi)$:
$$
\frac{\partial \psi}{\partial \nu} = \text{const.} \quad \text{on } F(\psi).
$$
In a seminal paper, Alt and Caffarelli \cite{AC} systematically analyzed this problem in the one-phase setting ($u \geq 0$) by minimizing the 
non-differentiable functional:
$$
J(v)=\int_{\Omega}|\nabla v|^2+\chi_{\{v>0\}}\,dx
$$
over an appropriate Sobolev class, proving that local minimizers satisfy 
$|\nabla u|^2 = 1$ along the free boundary $F(u) := \partial\{u>0\} \cap \Omega$. 
The two-phase framework, where solutions can change sign, was subsequently 
developed by Alt, Caffarelli, and Friedman \cite{ACF} using a groundbreaking monotonicity formula.
\begin{figure}
\begin{center}
\includegraphics[width=3in, scale=0.3]{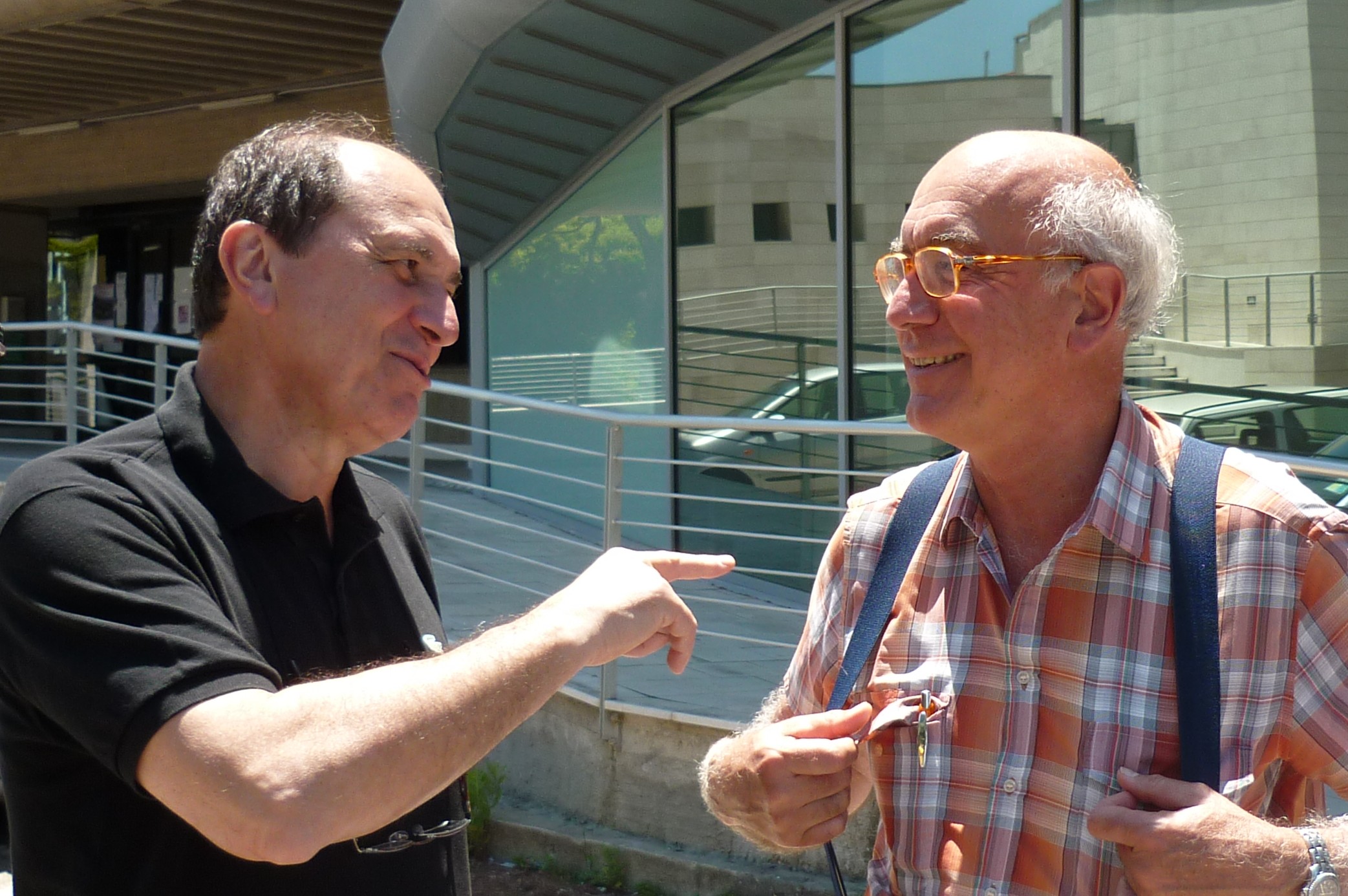}
\end{center}
\nopagebreak{
{\begin{center}Sandro Salsa (left) and Giorgio Talenti (right). Cagliari 2010  \end{center}}}
\end{figure}

Alongside the Stefan and Bernoulli problems, the obstacle problem serves as a core model case, describing the equilibrium configuration of an elastic membrane constrained to lie above a given profile $\psi$. Formally represented via variational inequalities, its core features can be illustrated by the localized equation $\Delta u = \chi_{\{u>0\}}$ with $u \geq 0$.
\begin{figure}
\begin{center}
\includegraphics[width=3in, scale=0.5]{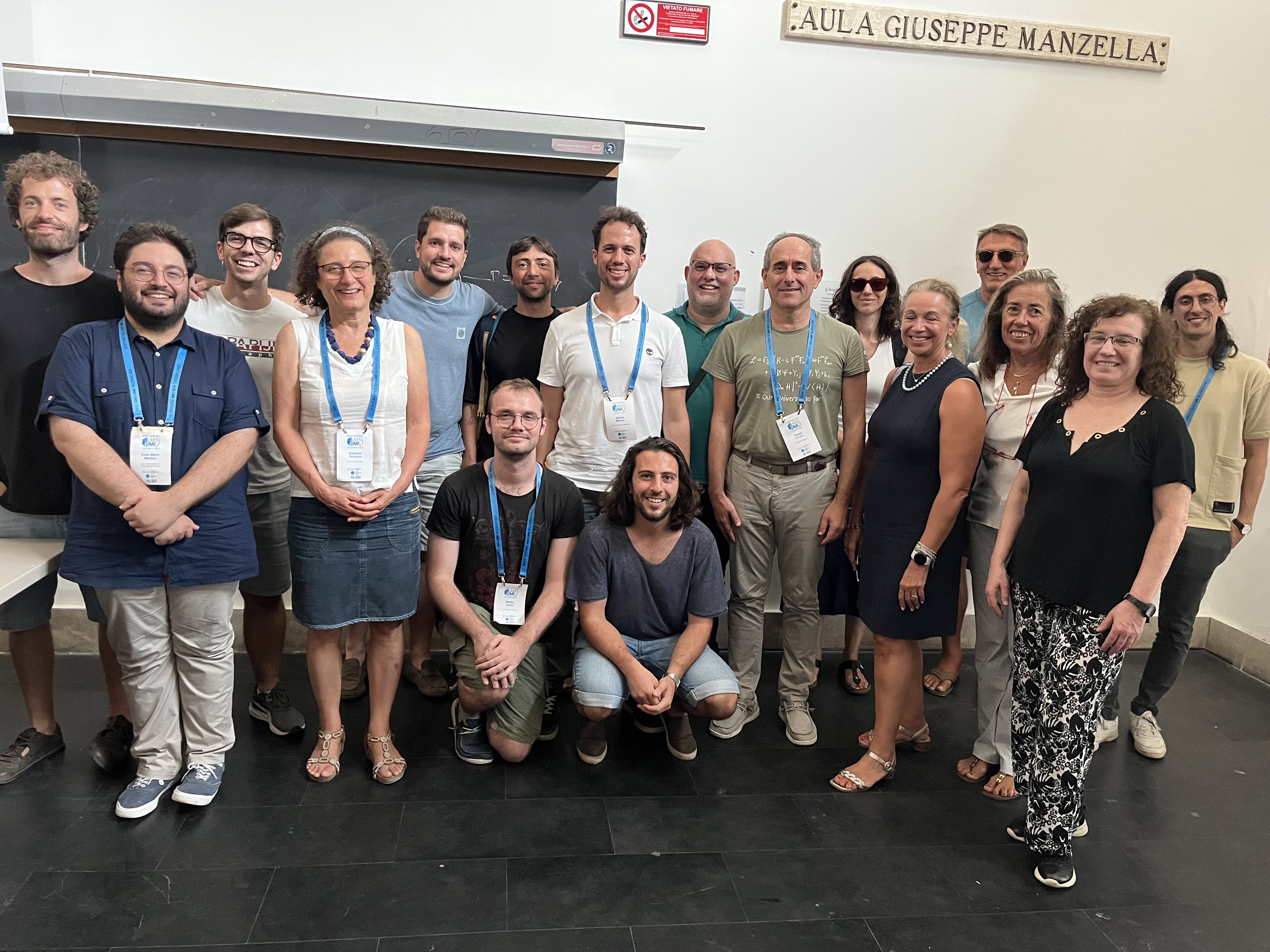}
\end{center}
\nopagebreak{
{\begin{center} AMS-UMI Meeting, Palermo, July 2024. Special Session on Free boundary problems: viscosity and variational approaches \end{center}}}
\end{figure}
The development of these three classical models deeply influenced the broader mathematical framework of the field. Throughout the 1950s to the 1970s, the variational and weak formulations for such problems were fully matured by P. Lax, H. Lewy, J. L. Lions, and G. Stampacchia, among others. This laid the groundwork for a major paradigm shift in the 1980s with the introduction of viscosity theory. Propelled by breakthrough results by Caffarelli \cite{C1, C2, C3}, a purely geometric approach was established to address the regularity of the free boundary under a viscosity formulation. On the parabolic side, this geometric regularity framework reached its peak with the fundamental trilogies \cite{ACS1, ACS2, ACS3} on the Stefan problem, and was later comprehensively synthesized in the landmark monograph \cite{CS}.
\begin{figure}
\begin{center}
\includegraphics[width=2in, scale=0.5]{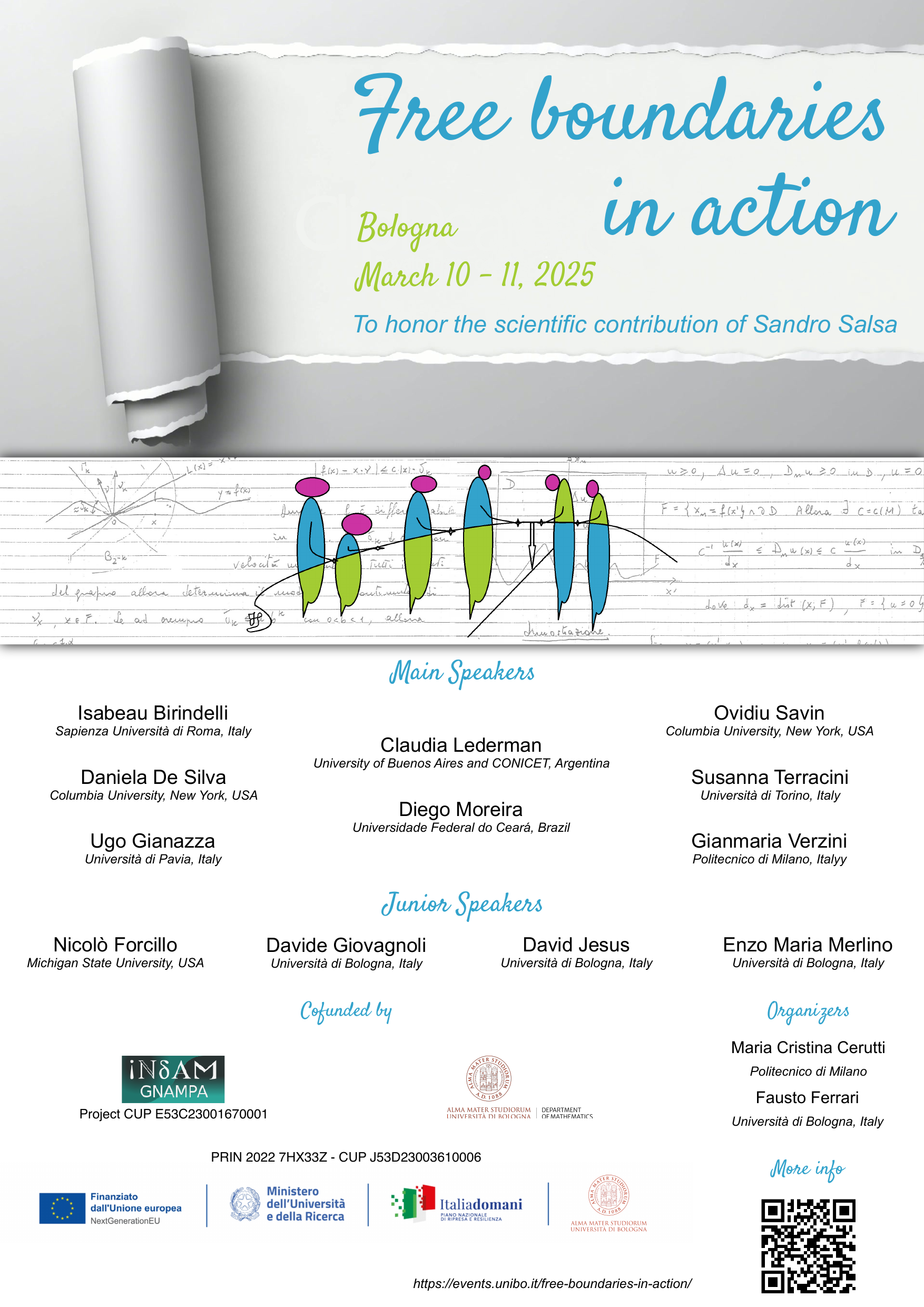}
\end{center}
\nopagebreak{
{\begin{center} Meeting in honor of the 75th birthday of Sandro Salsa. Bologna, March 2025  \end{center}}}
\end{figure}
In the decades since, these geometric techniques have expanded rapidly, see for instance \cite{CFS1, DS3, Fe, FeSa,F1,F2,FRS, LN1, LN2, LW,  W1, W2, W3, CPQ}. The classical linear theory has been extended to non-homogeneous one-phase problems \cite{D} and complex two-phase situations \cite{DFS, DFS1, DFS2, DFS3, MS, SaTV, DT}. Modern developments now encompass nonlinear operators such as the $p(x)$-Laplacian \cite{FeLe, FeLe2, FeLe3} and \cite{FeLe4} for the $p-$Laplace, multi-phase systems featuring double discontinuities \cite{DeSV}, thin free boundary variants restricted to lower-dimensional spaces \cite{CRS, DR, DS1, DS2, DSV,STV}, and sub-elliptic problems in non-commutative geometric structures \cite{BV, DGP, DGS, G}. Furthermore, the tools developed for the Bernoulli problem remain linked to active fields like the Hele-Shaw evolution equations \cite{CJK1, CJK2} and optimal partition problems \cite{KL, MTV}.

Sandro Salsa played a fundamental role in this development from its early stages. His foundational results laid the groundwork for the discipline's modern era and his own ongoing research, establishing the basic machinery that underpins the wide variety of free boundary problems studied today.

\vspace{1cm}

{\it Thank you Sandro!}


\begin{thebibliography}{30}
 

\bibitem{ACS1}I. Athanasopoulos, L. Caffarelli, S. Salsa,  \emph{Caloric functions in Lipschitz domains and the regularity of solutions to phase transition problems.} Ann. of Math. (2) 143 (1996), no. 3, 413--434.       
\bibitem{ACS2}I. Athanasopoulos, L. Caffarelli, S. Salsa,  \emph{Regularity of the free boundary in parabolic phase-transition problems.} Acta Math. 176 (1996), no. 2, 245--282.

\bibitem{ACS3}I. Athanasopoulos, L. Caffarelli, S. Salsa,  \emph{Phase transition problems of parabolic type: flat free boundaries are smooth.} Comm. Pure Appl. Math. 51 (1998), no. 1, 77--112.


\bibitem{AC} H.W. Alt, L.A.  Caffarelli,
\emph{Existence and regularity for a minimum problem with free
boundary.} J. Reine Angew. Math {\bf 325} (1981),105--144.

\bibitem{BV}I. Birindelli, E. Valdinoci,
\emph{The Ginzburg–Landau equation in the Heisenberg group.}
Commun. Contemp. Math. 10 (2008), no. 5, 671–719.

  \bibitem{C1} L. A.  Caffarelli, \emph{A Harnack inequality approach to the regularity of free boundaries. I. Lipschitz free boundaries are $C^{1,\alpha}$. }
Rev. Mat. Iberoamericana 3 (1987), no. 2, 139--162. 
  \bibitem{C2} L. A.  Caffarelli, \emph{A Harnack inequality approach to the regularity of free boundaries. II. Flat free boundaries are Lipschitz.} Comm. Pure Appl. Math. 42 (1989), no. 1, 55--78. 
 \bibitem{C3} L. A.  Caffarelli, \emph{A Harnack inequality approach to the regularity of free boundaries. III. Existence theory, compactness, and dependence on X.} Ann. Scuola Norm. Sup. Pisa Cl. Sci. (4) 15 (1988), no. 4, 583--602 (1989).

 \bibitem{CFMS} L. A.  Caffarelli, E. Fabes, S.  Mortola, S. Salsa,\emph{ Boundary behavior of nonnegative solutions of elliptic operators in divergence form.} Indiana Univ. Math. J. 30 (1981), no. 4, 621--640. 

\bibitem{ACF}  H.W. Alt, L.A.  Caffarelli,, A. Friedman,  \emph{Variational problems with two phases and their free boundaries. }Trans. Amer. Math. Soc. {\bf282} (1984), no. 2, 431--461.
\bibitem{CPQ} L. Caffarelli, S. Patrizi, V.  Quitalo, V. \emph{On a long range segregation model.} J. Eur. Math. Soc. (JEMS) 19 (2017), no. 12, 3575–3628.
   \bibitem{CRS} L.A. Caffarelli J-M. Roquejoffre, Y. Sire, \emph{ Variational problems with free boundaries for the fractional Laplacian}, J. European Math. Soc.,
{\bf 12} (2010), 1151--1179.

\bibitem{CS} L. Caffarelli, S. Salsa,  \emph{A geometric approach to free boundary problems.} Graduate Studies in Mathematics, 68. American Mathematical Society, Providence, RI, 2005.
\bibitem{CFS1} M. C. Cerutti, F. Ferrari, S. Salsa, \emph{Two-phase problems
for linear elliptic operators with variable coefficients: Lipschitz free
boundaries are $C^{1,\gamma }$.} Arch. Rational Mech. Anal. 171 (2004),
329-348

\bibitem{CSV}G-Q Chen, H. Shahgholian, J-L Vazquez, \emph{ Free boundary problems: the forefront of current and future developments},  Philosofical Transactions A 373 (2050), 2015


\bibitem{CJK1} S. Choi, D, Jerison, I. Kim, \emph{ Regularity for the one-phase hele-shaw problem from a lipschitz initial surface.} American journal of mathematics, 129(2):527--582, 2007.
\bibitem{CJK2} S, Choi, D. Jerison, I. Kim, \emph{ Local regularization of the one-phase hele-shaw flow.} Indiana University mathematics journal, pages 2765--2804, 2009.
\bibitem{DGP} D. Danielli, N. Garofalo, A. Petrosyan, \emph{The sub-elliptic obstacle problem: $C^{1,\alpha}$ regularity of the free boundary in Carnot groups of step two}. Adv. Math. 211 (2007), no. 2, 485--516.
\bibitem{DGS} D. Danielli, N. Garofalo, S. Salsa, \emph{Variational inequalities with lack of ellipticity. I. Optimal interior regularity and non-degeneracy of the free boundary}. Indiana Univ. Math. J. 52 (2003), no. 2, 361--398.
\bibitem{DeSV}  G. De Philippis, L. Spolaor, B. Velichkov, \emph{Regularity of the free boundary for the two-phase Bernoulli problem}. Invent. Math. 225 (2021), no. 2, 347--394.
\bibitem{D} D. De Silva, \emph{Free boundary regularity for a problem with
right hand side,} Interfaces Free Bound. 13 (2011), no. 2, 223--238.

\bibitem{DFS} D. De Silva, F. Ferrari, S. Salsa, \emph{Two-phase problems
with distributed sources: regularity of the free boundary,} Anal. PDE 7
(2014), no. 2, 267--310.

\bibitem{DFS1} D. De Silva, F. Ferrari, S. Salsa ,\emph{Free boundary
regularity for fully nonlinear non-homogeneous two-phase problems,} J. Math.
Pures Appl. (9) 103 (2015), no. 3, 658--694.

\bibitem{DFS2} D. De Silva, F. Ferrari, S. Salsa, \emph{Perron's solutions
for two-phase free boundary problems with distributed sources,} Nonlinear
Anal. 121 (2015), 382--402.

\bibitem{DFS3} D. De Silva, F. Ferrari, S. Salsa, \emph{Regularity of
higher order in two-phase free boundary problems,} Trans. Amer. Math. Soc.
 Trans. Amer. Math. Soc. 371 (2019), no. 5, 3691--3720




\bibitem{DR} De Silva D., Roquejoffre J-M., \emph{ Regularity in a one-phase  free boundary problem for the fractional Laplacian.}  Ann. Inst. H. Poincare Anal. Non Lineaire 29 (2012), no. 3, 335--367. 

\bibitem{DS1} De Silva D., Savin O., \emph{$C^{2,\alpha}$ regularity of flat free boundaries for the thin one-phase problem,}  Journal of Diff. Equations 253 (2012), no. 8, 2420--2459.

\bibitem{DS2} De Silva D., Savin O., \emph{Regularity of Lipschitz free boundaries for the thin one-phase problem},  J. Eur. Math. Soc. (JEMS) 17 (2015), no. 6, 1293--1326.

\bibitem{DS3} D. De Silva, O. Savin, \textit{Lipschitz regularity of
solutions to two-phase free boundary problems}, Int. Math. Res. Notices (2019), no. 7, 2204--2222.

\bibitem{DSV} S. Dipierro, O. Savin, E. Valdinoci, \emph{A nonlocal free boundary problem} SIAM Journal on Mathematical Analysis 47 (2015), no. 6, 4559--4605.

\bibitem{DT} D. De Silva, G. Tortone, \emph{Improvement of flatness for vector valued free boundary problems.} Math. Eng. 2 (2020), no. 4, 598–613. 


\bibitem{FGS1} E. B. Fabes, N. Garofalo, S. Salsa,   \emph{ Comparison theorems for temperatures in noncylindrical domains}. Atti Accad. Naz. Lincei Rend. Cl. Sci. Fis. Mat. Nat. (8) 77 (1984), no. 1--2, 1--12 (1985)
\bibitem{FGS2} E. B. Fabes, N. Garofalo, S. Salsa,   \emph{A backward Harnack inequality and Fatou theorem for nonnegative solutions of parabolic equations}. Illinois J. Math. 30 (1986), no. 4, 536--565.

\bibitem{FGMS} E. B. Fabes, N. Garofalo, S. Marine-Malave, S. Salsa,   \emph{Fatou theorems for some nonlinear elliptic equations}. Rev. Mat. Iberoamericana 4 (1988), no. 2, 227--251. 

\bibitem{Fe} F. Ferrari, \emph{Two-phase problems for a class of fully
nonlinear elliptic operators. Lipschitz free boundaries are $C^{1,\gamma }$}%
, American Journal of Mathematics, 128 (2006) no.3, 541-571.
 \bibitem{FeLe} F. Ferrari, C. Lederman,  \emph{Regularity of flat free boundaries for a p(x)-Laplacian problem with right hand side}. Nonlinear Anal. 212 (2021), Paper No. 112444, 25 pp.

\bibitem{FeLe2} F. Ferrari, C. Lederman, \emph{Regularity of Lipschitz free boundaries for a p(x)-Laplacian problem with right hand side}. J. Math. Pures Appl. (9) 171 (2023), 26--74
 \bibitem{FeLe3} F. Ferrari, C. Lederman, \emph{Regularity of flat free boundaries for two-phase $p(x)$-Laplacian problems with right hand side}, Calc. Var. Partial Differential Equations 63 (2024), no. 5, Paper No. 132, 43 pp.
 \bibitem{FeLe4} F. Ferrari, C. Lederman, \emph{Lipschitz regularity of solutions to two-phase $p$-Laplacian free boundary problems with right hand side}. Comm. Anal. Geom. 34 (2026), no. 2, 465--503
\bibitem{FeSa} F. Ferrari, S. Salsa, \emph{Regularity of the free boundary
in two-phase problems for linear elliptic operators}, Advances in Math, 214
(2007), 288-322.

\bibitem{F1} M. Feldman, \emph{Regularity for nonisotropic two-phase
problems with Lipschitz free boundaries.} Differential Integral Equations 10
(1997), no. 6, 1171--1179.

\bibitem{F2} M. Feldman, \emph{Regularity of Lipschitz free boundaries in
two-phase problems for fully nonlinear elliptic equations}, Indiana Univ.
Math. J. 50 (2001), no. 3, 1171--1200.
\bibitem{FRS} A. Figalli, X. Ros-Oton, J. Serra,  \emph{The singular set in the Stefan problem.} J. Amer. Math. Soc. 37 (2024), no. 2, 305--389. 
%
%

\bibitem{Friedman}A.  Friedman, \emph{ Variational principles and free-boundary problems}. A Wiley-Interscience Publication. Pure and Applied Mathematics. John Wiley\& Sons, Inc., New York, 1982.

\bibitem{G} U. Gianazza, \emph{Local properties of variational solutions for the two obstacle problem involving square Hörmander operators.} Ann. Mat. Pura Appl. 164 (1993), 301--333.

\bibitem{HW} R. A. Hunt, R.L. Wheeden, 
\emph{ Positive harmonic functions on Lipschitz domains.}
Trans. Amer. Math. Soc. 147 (1970), 507--527.

\bibitem{JK} D. Jerison, C. Kenig, \emph{ Boundary behavior of harmonic functions in nontangentially accessible domains.} Adv. in Math. 46 (1982), no. 1, 80--147.
\bibitem{Kam}
S. L. Kamenomostskaya,
\emph{On Stefan problem.}
Nauchnye Doklady Vysshey Shkoly, Fiziko-Matematicheskie Nauki
\textbf{1}, (1958) 60--62. [In Russian]

\bibitem{Kind-Stam-book}
D. Kinderlehrer,  G. Stampacchia,
\textit{An Introduction to Variational Inequalities and Their Applications}.
Pure and Applied Mathematics, \textbf{88}. Academic Press, Inc. [Harcourt Brace Jovanovich, Publishers]:
New York-London (1980).
\bibitem{KL} D. Kriventsov, F. Lin, \emph{Regularity for Shape Optimizers:
The Nondegenerate Case} CPAM  Comm. Pure Appl. Math. 71 (2018), no. 8, 1535--1596
\bibitem{LN1} J. Lewis, K.  Nystr\"om,\emph{ Regularity of flat free boundaries in two-phase problems for the p-Laplace operator}, Ann. Inst. H. Poincar\'e Anal. Non Lin\'aire 29 (2012), no. 1, 83--108.
\bibitem{LN2} J. Lewis, K.  Nystr\"om, \emph{Regularity and free boundary regularity for the p-Laplace operator in Reifenberg flat and Ahlfors regular domains}. J. Amer. Math. Soc. 25 (2012), no. 3, 827--862.
\bibitem{LW} Lederman C., Wolanski N., \emph{A two phase elliptic singular
perturbation problem with a forcing term,} J. Math. Pures Appl. (9) 
86 (2006), no. 6, 552--589.
\bibitem{LQT} R. Leitao, O. de Queiroz, E. Teixeira,
\emph{Regularity for degenerate two-phase free boundary problems}. 
Ann. Inst. H. Poincar\'e Anal. Non Lin\' eaire 32 (2015), no. 4, 741--762. 

\bibitem{MTV} D. Mazzoleni, S. Terracini, B. Velichkov,  \emph{ Regularity of the optimal sets for some spectral functionals}, Geom. Funct. Anal. 27, 373–426 (2017).

\bibitem{MS} D. Moreira, H. Shrivastava, \emph{Optimal regularity for variational solutions of free transmission problems.} Journal de Mathématiques Pures et Appliquées 169 (2023), 1–49.

\bibitem{SaTV} S. Salsa, F. Tulone, G. Verzini,  \emph{Existence of viscosity solutions to two-phase problems for fully nonlinear equations with distributed sources.} Math. Eng. 1 (2019), no. 1, 147-173

\bibitem{STV} Y. Sire, S. Terracini, S. Vita, \emph{Liouville type theorems and regularity of solutions to degenerate or singular problems.Part I: Even solutions.} Communications in Partial Differential Equations 46 (2021), no. 2, 310–361. 


\bibitem{S} J. Stefan, \emph{\"Uber die Theorie der Eisbildung, insbesondere über die Eisbildung im Polarmeere}, Sitzungsberichte der \"Osterreichischen Akademie der Wissenschaften, Mathematisch-Naturwissenschaftliche Klasse,
Abteilung 2, Mathematik, Astronomie, Physik, Meteorologie und Technik,(1898) 98, 965--983



\bibitem{W1} P.Y. Wang, \emph{Regularity of free boundaries of two-phase
problems for fully nonlinear elliptic equations of second order. I.
Lipschitz free boundaries are $C^{1,\alpha }$}, Comm. Pure Appl. Math. 53
(2000), no. 7, 799--810.

\bibitem{W2} P.Y. Wang, \emph{Regularity of free boundaries of two-phase
problems for fully nonlinear elliptic equations of second order. II. Flat
free boundaries are Lipschitz}, Comm. in Partial Differential equations
27(7\&8) (2002), no. 7, 1497--1514.

\bibitem{W3} P.Y. Wang, \emph{Existence of solutions of two-phase free
boundary for fully non linear equations of second order, }J. of Geometric
Analysis (2002), no. 7, 1497--1514.


\end{thebibliography}
\end{document}